\documentclass{ifacconf}

\usepackage{graphicx}      
\usepackage{natbib}        

\usepackage{amsmath, amssymb}
\usepackage{algorithm}
\usepackage{algpseudocode}
 \usepackage{multirow}

\newcommand*{\R}{{\mathbb R}}
\def\x{\mathbf{x}}
\def\y{\mathbf y}

\DeclareMathOperator*{\argmin}{argmin}
\DeclareMathOperator*{\argmax}{argmax}

\usepackage{xcolor}

\begin{document}
\begin{frontmatter}

\title{Adaptive Mirror Descent for the Network Utility Maximization Problem\thanksref{footnoteinfo}} 

\thanks[footnoteinfo]{The research of F. Stonyakin in Section 3 was supported by Russian Science Foundation project 18-71-00048. The research of E. Vorontsova was supported by Russian Foundation for Basic Research project 18-29-03071 mk. The research of A. Ivanova, F. Stonyakin, A. Gasnikov in Sections 4 and 5 was supported by Russian Foundation for Basic Research project 19-31-51001.
The research of A. Gasnikov was partially supported by Yahoo! Research Faculty Engagement Program.}

\author[First]{Anastasiya Ivanova} 
\author[Second]{Fedor Stonyakin} 
\author[Third]{Dmitry Pasechnyuk}
\author[Fourth]{Evgeniya Vorontsova}
\author[Fifth]{Alexander Gasnikov}

\address[First]{Moscow Institute of Physics and Technology, Moscow, Russia, \\ National Research University Higher School of Economics, Moscow, Russia, Sirius University of Science and Technology, Sochi, Russia (anastasiya.s.ivanova@phystech.edu)
} 

\address[Second]{Moscow Institute of Physics and Technology, Moscow, Russia,\\ 
V. I. Vernadsky Crimean Federal University, Simferopol, Russia, \\ Sirius University of Science and Technology, Sochi, Russia}
\address[Third]{Saint Petersburg Lyceum 239, Russia}
\address[Fourth]{Grenoble Alpes University, Grenoble, France}
\address[Fifth]{Moscow Institute of Physics and Technology, Moscow, Russia,\\
Institute for Information Transmission Problems, Moscow, Russia,\\
Caucasus Mathematical Center, Adyghe State University, Russia, Sirius University of Science and Technology, Sochi, Russia}

\begin{abstract}   
Network utility maximization is the most important
problem in network traffic management.
 Given the
growth of modern 
communication networks,
we consider utility
maximization problem in a network with a large number of connections (links) that are used by a huge number of users. 
To solve this problem an adaptive mirror descent 
algorithm for many constraints is proposed. 
The key feature of the algorithm
is that it has a dimension-free convergence rate. 
The convergence of the proposed scheme is
proved theoretically.
The theoretical analysis is verified with numerical simulations.
We compare the algorithm with another
approach, using the ellipsoid method~(EM) for the dual problem.
Numerical experiments showed that
the performance of the proposed algorithm against EM
is significantly better in large networks and when very high solution accuracy
is not required.
Our approach can be used
in many network design paradigms, in particular,
in software-defined networks.
\end{abstract}


\end{frontmatter}

\section{Introduction}
One of the most significant design features of modern 
communication
network
systems is the capacity  
to adjust the distribution of bandwidth
and other network properties
for achieving best performance 
and reliability in real-time.
The tasks of finding the best configuration
and/or design parameters for networks
are actually reduced to solving 
complex optimization problems with thousands and
millions of variables. 

The key question we address in this article is how the available bandwidth on the network should be distributed among competing connections. In this case, users can control the use of available bandwidth by adjusting the connection price. 
As mentioned earlier, 
in view of the huge size of such networks it is very important to develop algorithms dimension-free
on
the size of the system. 

Thus, network utility maximization~(NUM) problems in computer networks with a large number of connections are
considered.
Connections are used for their own purposes
by consumers (users), the number of whom
can also be very large. 
The purpose of the work is to determine the mechanism of resource allocation, which in the context of this task are available bandwidth connections. At the same time, it is necessary to ensure stable operation of the system and prevent overloads. As an optimality criterion, the sum of the utilities of all users of the network is used.

The original resource allocation framework,
reduced to the maximization
of aggregate concave utility functions
subject to link capacity constraints, was pioneered by~\cite{kelly}.
During the last two decades 
NUM framework has found wide-ranging applications
to wireless and sensor networks, and many
other fields~\cite{palomar, dehghan},
for a survey, see~\cite{shak2008survey} and references therein.

Also, the 
mechanisms of decentralized resource allocation proposed in a monograph by~\cite{arrow1958decentralization} have since attracted much attention in economic research, see, for example, \cite{kakhbod2013resource, campbell1987resource, friedman1995complexity}. 
Furthermore, the problem of resource distribution
in computer networks 
was investigated in a recent paper by~\cite{rokhlin} as well. In this paper, following~\cite{nesterov2018dual, ivanova2019res, ivanova2018composite}, we additionally consider a price adjustment mechanism. The proposed approach has a practical value due to  decentralization, which means
that to set and adjust the price of an individual connection, only the reaction of the users who use that connection is necessary,
not the reaction of all users. 

\section{Problem Statement}
Consider a communication network with $m$~connections (links)
and $n$~users (or nodes).
Users
exchange packets through a fixed 
connections' set. 
The network structure is given by the binary routing matrix $C = (C_i^j)\, \in \, \R^{m \times n}$.
The columns of the matrix $\mathbf{C}_i \neq 0$, $i = 1,\, \ldots,
\, n$ are boolean $m$-dimensional vectors such that
$C_i^j = 1$ if 
node~$i$ is used in the connection~$j$, otherwise $C_i^j = 0$.
Capacity constraints are given by the vector~$\mathbf{b} \, \in \, \R^m$
with strictly positive components.
These constraints imply that no connection is
overloaded.
Users evaluate the quality of the network
with utility functions $u_k (x_k)$, $k = 1,\, \ldots,\, n$, where $x_k \, \in \, \R_+$ is the network
rate for the $k$--th user.

The problem of maximizing the total utility of the network under the given constraints is formulated as follows:
\begin{eqnarray}
\label{eq_main_task}
    \max_{\left \{ C \mathbf{x} = \sum\limits_{k = 1}^n \mathbf{C}_k x_k \right \} \, \le \, \mathbf{b}} \left \{ U( \mathbf{x}) = \sum_{k = 1}^n u_k(x_k) \right \},
\end{eqnarray}
where $u_k(x_k)\, \, \, k = 1, \ldots, n$ are concave functions and $\mathbf{x} = (x_1, \, \ldots, \, x_n) \, \in \, \mathbb{R}^n_+$. 
Optimal
resource allocation~$\mathbf{x^*}$ is a solution
of the problem~\eqref{eq_main_task}.

Denote $g_j(\x) = \langle \mathbf{C}_{j},\x \rangle - b_{j},\, \, \, j = 1, \ldots, m$. Note that 
\begin{eqnarray*}
    |g_j(\x^1) - g_j(\x^2)| \leq \|\mathbf{C}_{j}\|_{2} \|\x^1 - \x^2\|_2 \leq m \cdot \|\x^1 - \x^2\|_2.
\end{eqnarray*}
The last inequality holds due to the definition of matrix~$C$. So, each $g_j(\x),\, \, \, j = 1, \ldots, m$ is $M_{g} = m$-Lip\-schitz con\-tinuous.

For convenience, we will further consider the minimization problem equivalent to~\eqref{eq_main_task} 
\begin{eqnarray}
\label{eq_main_task_min}
    \min_{ g_j(\x) \leq 0, \, j = \overline{1, \, m} } f(\x),
\end{eqnarray}
where $f(\x) = - U(\x)$.

\section{Mirror Descent for Many Constraints}\label{sect_mir}
Let us consider the max-type functional constraint: $g(x)=\max\limits_{j\in \overline{1,m}}g_j(x),$ which keeps the Lipschitz property and non-smoothness, provided that all functions $g_j(x)$,
$j = 1, \, \ldots, \, m$, satisfy these properties. 

\subsection{First algorithm}
\label{Mirr_Des_1}
The first variant of Mirror Descent Algorithm 
for many constraints is suggested,
see Algorithm~\ref{alg4m}.

\begin{algorithm}
\caption{Mirror Descent Algorithm for many constraints}
\label{alg4m}
\begin{algorithmic}[1]
\Require $\varepsilon>0,\Theta_0: \,d(\x^*) = \frac{1}{2}\|\x^0 - \x^* \|_2 \leqslant\Theta_0^2$, initial point~$\x^0 = 0$
\State $I =: \emptyset$
\State $N \leftarrow 0$
\Repeat
\If
{$g_j(\x^N)\leqslant\varepsilon\|\nabla g_{j}(\x^N)\|_2, \, \forall \, j = \overline{1, \, m}$}
\State $\x^{N+1}=\left[\x^N - \frac{\varepsilon \nabla f(\x^N)}{\|\nabla f(\x^N)\|_2} \right]_{+}$ 
\State // $h_N=\frac{\varepsilon}{\|\nabla f(x^N)\|_2}$
\State $N \rightarrow I$
\Else
\State \emph{$(g_{j_N}(\x^N)>\varepsilon\|\nabla g_{j_N}(\x^N)\|_2), \, \, j_N = \overline{1, \, m})$}
\State $\x^{N+1}=\left[\x^N - \frac{\varepsilon \nabla g_{j_N}(\x^N)}{\|\nabla g_{j_N}(\x^N)\|_2} \right]_{+}$ 
\State // $h_N=\frac{\varepsilon}{\|\nabla g_{j_N}(\x^N)\|_2}$
\EndIf
\State $N\leftarrow N+1$

\Until {$2\frac{\Theta_0^2}{\varepsilon^2}\leqslant N$}
\Ensure $\bar{\x}^N:= \argmin\limits_{\x^k,\;k\in I}\,f(\x^k)$
\end{algorithmic}
\end{algorithm}

Now we will estimate the rate of convergence of the proposed method. For this we need the following auxiliary assumption (\cite{nesterov2018lectures}, Lemma 3.2.1). Recall that $\x^*$~ is the solution of the problem~\eqref{eq_main_task_min}.


\begin{lem}
\label{LemmT}
Let us define the following function:
\begin{equation} \label {eq13}
\omega (\tau) = \max \limits_ {\x \in \mathbb{R}^n_+} \{f (\x) -f (\x^*): \| \x-\x^* \| \leqslant \tau \},
\end{equation} where $\tau$ is a positive number.
Then for any $\y \in X$
\begin{equation} \label {eq_lemma}
f (\y) - f (\x^*) \leqslant \omega (v_f (\y, \x^*)),
\end{equation}
where
\begin{eqnarray*}
    v_{f}(\y, \x^*) = \left \langle \frac{\nabla f(\y)}{\|\nabla f(\y) \|}, \y - \x^* \right \rangle  \, \, \text{for} \, \nabla f(\y) \neq 0
\end{eqnarray*}
and $v_{f}(\y, \x^*) = 0$ for $\nabla f(\y) = 0$.
\end{lem}

\begin{thm}\label{th2_new_methods}
Let $\varepsilon > 0$ be a fixed number and the stopping criterion of Algorithm \ref{alg4m} be satisfied. Then
\begin{equation}\label{eq1_new_methods}
\min\limits_{k \in I} v_f(\x^k,\x^*) \leqslant \varepsilon, \ \ \max\limits_{k\in I} g_j(\x^k)\leq\varepsilon M_g, \, \, j\in \overline{1,m}.
\end{equation}
\end{thm}

\begin{pf}

1) If $k\in I$ (for productive steps),
\begin{eqnarray}
\label{eq2}
& & h_k(f(\x^k)-f(\x^*)) \leq h_k  \langle \nabla f(\x^k),\x^k-\x^* \rangle = \varepsilon v_f(\x^k,\x^*) \notag \\
& \leqslant & \frac{h_k^2}{2}||\nabla f(\x^k)||_2^2+\frac{1}{2}\|\x^k - \x^* \|^2_2 -\frac{1}{2}\|\x^{k+1} - \x^* \|^2_2 \notag\\ &
=&\frac{\varepsilon^2}{2}+\frac{1}{2}\|\x^k - \x^* \|^2_2 -\frac{1}{2}\|\x^{k+1} - \x^* \|^2_2.
\end{eqnarray}

2) If $k\not\in I$, then $\frac{g_{j_k}(\x^k)-g_{j_k}(\x^*)}{||\nabla g_{j_k}(\x^k)||_2} \geq \frac{g_{j_k}(\x^k)}{||\nabla g_{j_k}(\x^k)||_2}>\varepsilon$. Therefore, the following inequalities hold
\begin{eqnarray}
\label{eq3}
&\varepsilon^2&<h_k(g_{j_k}(\x^k)-g_{j_k}(\x^*))\leqslant\frac{h_k^2}{2}||\nabla g_{j_k}(\x^k)||_2^2 \notag \\ &+& \frac{1}{2}\|\x^k - \x^* \|^2_2  \notag  -\frac{1}{2}\|\x^{k+1} - \x^* \|^2_2 \notag \\ &=&\frac{\varepsilon^2}{2}+\frac{1}{2}\|\x^k - \x^* \|^2_2-\frac{1}{2}\|\x^{k+1} - \x^* \|^2_2, \notag \\ &\text{ or}&
 \quad \frac{\varepsilon^2}{2}<\frac{1}{2}\|\x^k - \x^* \|^2_2-\frac{1}{2}\|\x^{k+1} - \x^* \|^2_2.
\end{eqnarray}

3) After summing the inequalities \eqref{eq2} and \eqref{eq3} we have:
\begin{eqnarray*}
&& \sum_{k\in I} \varepsilon v_f(\x^k,\x^*) \leqslant |I|\frac{\varepsilon^2}{2}-\frac{\varepsilon^2|J|}{2} \\ &+& \frac{1}{2}\|\x^0 - \x^* \|^2_2- \frac{1}{2}\|\x^{k+1} - \x^* \|^2_2
=\varepsilon^2|I|-\frac{\varepsilon^2 N}{2}+\Theta_0^2.
\end{eqnarray*}
And since $\sum_{k\in I} v_f(\x^k,\x^*)  \geq |I|\min\limits_{k \in I} v_f(\x^k,\x^*)$, after the stopping criterion of the algorithm holds we have
\begin{eqnarray*}
|I| \min\limits_{k \in I} \varepsilon v_f(\x^k,\x^*) &\leqslant& \varepsilon^2|I|-\frac{\varepsilon^2 N}{2}+\Theta_0^2  \leq \varepsilon^2|I|. 
\end{eqnarray*}
So,  $\min\limits_{k \in I} v_f(\x^k,\x^*) \leqslant \varepsilon.$

Further, for each $k\in I\;g_j(\x^k)\leqslant\varepsilon||\nabla g_j(\x^k)||_2 \leqslant\varepsilon M_g, \,  j = 1, \ldots, m$.

Now we have to show that the set of productive steps $I$ is non-empty. If $I=\emptyset$, then $|J|=N$ and  the Lipschitz continuous of $g$ means that $N\geqslant\frac{2\Theta_0^2}{\varepsilon^2}$. On the other hand, from \eqref{eq3} we have:
$$\frac{\varepsilon^2N}{2}< \frac{1}{2}\|\x^0 - \x^* \|^2_2\leqslant\Theta_0^2,$$
which leads us to the controversy, so $I\neq\emptyset$.
\end{pf}
Now let us show how to estimate the quality of the solution by the function basing on the previous theorem for Lipschitz continuous function. 

\begin{cor}
\label{cor_lipschits1}
Let $f$ satisfy the Lipschitz condition
\begin{equation}\label{lipschits_condition}
|f(\x)-f(\y)|\leqslant M_f\|\x-\y\|_2 \quad \forall \x,\y\in X.
\end{equation}
Then, after the stopping of Algorithm \ref{alg4m}, the following inequality holds:
$$\min\limits_{k\in I}f(\x^k)-f(\x^*)\leq M_f \varepsilon.$$
\end{cor}
\begin{pf}
Note that 
\begin{eqnarray*}
\min\limits_{k\in I}f(\x^k)-f(\x^*) \leq \min\limits_{k \in I} v_f(\x^k,\x^*) \cdot \|\nabla f(\x^k) \|_2 \leq M_{f}\varepsilon.
\end{eqnarray*}
\end{pf}

Now, we estimate the rate of convergence of Algorithm~\ref{alg4m} for a differentiable objective functional $f$ with a Lipschitz-continuous gradient.
\begin {equation} \label{eqlipgrad}
\|\nabla f (x) - \nabla f (y)\|_* \leqslant L \|x-y\| \quad \forall x, y \in X.
\end {equation}
Assume that similarly to \cite{devolder2014first} we have an inexact $(\delta,L)$-gradient 
$\nabla_{\delta} f$ for $f$:
$$ f (\x) \leqslant f (\x_*) + \langle \nabla_{\delta} f (\x_*), \x-\x_* \rangle + \frac {1}{2} L \| \x-\x_* \|^2 + \delta $$
for exact solution $\x_*$ we can get that
\begin{eqnarray*}
&& \min\limits_ {k \in I} f (\x^k) -f (\x_*)\\ & \leqslant &  \min \limits_ {k \in I} \left \{\| \nabla_{\delta} f (\x_*) \|_2 \| \x^k-\x_* \| + \frac{1}{2} L \| \x^k-\x_* \| ^ 2  + \delta  \right \}.
\end{eqnarray*}
From  Theorem~\ref{th2_new_methods} in view of Lemma~\ref{LemmT} we  obtain the following corollary 
\begin{cor}
\label{cor_lipschits_grad}
Let $f$ be differentiable and~\ref{eqlipgrad} hold. Assume that we have an inexact $(\delta,L)$-gradient $\nabla_{\delta} f$ of function $f$ at each point x. Then, after the stopping of Algorithm \ref{alg4m}, the next inequality holds:
\begin{eqnarray*}
 \min\limits_ {k \in I} f (\x^k) - f (\x_*) \leqslant \varepsilon \| \nabla f (\x_*) \|_2  + \frac{1}{2} L \varepsilon^2  + \delta .
\end{eqnarray*}
\end{cor}

Let us consider the rate of convergence of Algorithm~\ref{alg4m} for a differentiable objective H\"{o}lder-continuous functional $f$, i.e. for some $\nu\in[0;1)$
\begin{equation}\label{gelder_condition}
|f(x)-f(y)|\leqslant M_{f,\;\nu}\|x-y\|^{\nu}\quad \forall x,y\in Q.
\end{equation}

For example, $\nu = 1/2$ for $f(x)=\sqrt{x}$. Let us recall the following inequality \cite{stonyakin2019universal}
\begin{equation}\label{eq002}
M_{\nu}a^{\nu}\leqslant M_{\nu}\left[\frac{M_{\nu}}{\delta}\right]^{\frac{1-\nu}{1+\nu}}\frac{a^2}{2}+\delta,
\end{equation}
which is true for each $\delta>0$. Then by \eqref{gelder_condition} we have
$$|f(x)-f(y)|\leqslant\frac{M_{\nu}^{\frac{2}{1+\nu}}}{2\delta^{\frac{1-\nu}{1+\nu}}}||x-y||_2^2+\delta.$$

Set $\delta=\varepsilon$. Then
\begin{equation}\label{eq003}
|f(x)-f(y)|\leqslant\underbrace{\frac{M_{\nu}^{\frac{2}{1+\nu}}}{2\varepsilon^{\frac{1-\nu}{1+\nu}}}}_{M}||x-y||_2^2+\varepsilon.
\end{equation}
By Lemma ~\ref{LemmT} after the stopping of Algorithm \ref{alg4m} we have $\min\limits_{k\in I}v_f(x^k,x_*)<\varepsilon$. It means the following inequality:
\begin{equation}\label{eq004}
f(\widehat{x})-f^*\leqslant\frac{M_{\nu}^{\frac{2}{1+\nu}}}{2\varepsilon^{\frac{1-\nu}{1+\nu}}}\varepsilon^2+\varepsilon =
\frac{M_{\nu}^{\frac{2}{1+\nu}}}{2} \varepsilon^{1+\frac{2\nu}{1+\nu}} + \varepsilon.
\end{equation}
Then we can formulate the following corollary 
\begin{cor}\label{CorHolder}
Let $f$ be a H\"{o}lder-continuous functional and~\eqref{gelder_condition} hold. Then, after the stopping of Algorithm \ref{alg4m} for $\varepsilon<1$ the inequality \eqref{eq004} means
$$f(\widehat{x})-f^* \leqslant \widehat{M}\varepsilon$$
for some $\hat{M} > 0$.
\end{cor}
So, for problems with a convex H\"{o}lder-continuous differentiable objective and convex Lipcshitz-continuous functional constraints we can achieve an $\varepsilon$-solution after
$
O\left(\frac{1}{\varepsilon^2}\right)
$
iterations of  Algorithm \ref{alg4m}. This estimate is optimal due to its optimality on a significantly narrower class of problems with Lipschitz-continuous objective functionals~\cite{nemirovsky1983problem}.

\subsection{Second algorithm}

Let us consider the following method for fixed accuracy $\varepsilon>0$, initial approach $x^0,\; \Theta_0$: $\|x^0 - x_*\|_2^2 \leqslant 2\Theta_0^2$ and Lipschitz-continuous functional constraint $g$:
$$|g(x)-g(y)|\leqslant M_g||x-y||\;\forall x,y\in \mathbb{R}^n_+.$$
\begin{algorithm}
\caption{Another variant of adaptive mirror descent for many constraints.}
\label{alg1}
\begin{algorithmic}[1]
\Require $\varepsilon>0,\Theta_0: \,d(\x^*) = \frac{1}{2}\|\x^0 - \x^* \|_2 \leqslant\Theta_0^2$, $\x^0=0$ - initial point.
\State $I=:\emptyset$
\State $N\leftarrow0$
\Repeat
\If
{$g_j(\x^N)\leqslant\varepsilon\|\nabla g_{j}(\x^N)\|_2, \, \forall j\in \overline{1,m}$}
\State $\x^{N+1}=\left[\x^N - \frac{\varepsilon \nabla f(\x^N)}{\|\nabla f(\x^N)\|^2_2} \right]_{+}$ 
\State // $h_N=\frac{\varepsilon}{\|\nabla f(x^N)\|^2_2}$
\State $N\rightarrow I$
\Else
\State \emph{$(g_{j_N}(\x^N)>\varepsilon\|\nabla g_{j_N}(\x^N)\|_2), \, \, j_N \in \overline{1,m})$}
\State $\x^{N+1}=\left[\x^N - \frac{\varepsilon \nabla g_{j_N}(\x^N)}{\|\nabla g_{j_N}(\x^N)\|_2} \right]_{+}$ 
\State // $h_N=\frac{\varepsilon}{\|\nabla g_{j_N}(\x^N)\|_2}$
\EndIf
\State $N\leftarrow N+1$

\Until {\begin{equation}\label{eq11}
\frac{2\Theta_0^2}{\varepsilon^2} \leqslant\sum\limits_{k\in I}\frac{1}{||\nabla f(x^k)||_{*}^2}+|J|,
\end{equation}

where $|J|$~--- the number of unproductive steps (we denote by $|I|$ the number of productive steps, i.e. $|I|+|J|=N$).}
\Ensure $\hat{\x}^N=\frac{1}{\sum_{k\in I}h_k}\sum\limits_{k\in I}h_k\x^k$
\end{algorithmic}
\end{algorithm}

\begin{thm}\label{th1}
Let $\varepsilon > 0$ be a fixed number and the stopping criterion of Algorithm~\ref{alg1} be satisfied. Then the following inequality is true:
\begin{eqnarray*}
f(\hat{\x}^N)-f(\x_*)\leqslant\varepsilon\text{ and} \label{th1_1} \\
g(\hat{\x}^N)\leqslant\frac{\varepsilon}{\sum_{k\in I}h_k}\sum_{k\in I}h_k||g(\x^k)||_2\leqslant\varepsilon M_g \label{th1_2} ,
\end{eqnarray*}
where $\hat{\x}^N=\frac{1}{\sum_{k\in I}h_k}\sum\limits_{k\in I}h_k\x^k$.
\end{thm}
\begin{pf}
We give only a sketch of the proof because the proof of this theorem mostly follows the proof of theorem~\ref{th2_new_methods}. 
1) If $k\in I$ (for productive steps),
\begin{eqnarray*}
& & h_k(f(\x^k)-f(\x^*))  \\ &
\leq &\frac{\varepsilon^2}{2}\cdot\frac{1}{||\nabla f(\x^k)||_2^2}+ \frac{1}{2}\|\x^k - \x^* \|^2_2- \frac{1}{2}\|\x^{k+1} - \x^* \|^2_2.
\end{eqnarray*}
2) If $k \not\in I$ inequality~\eqref{eq3} holds.

3) Summing the inequalities and after fulfilling the criterion for stopping the algorithm \eqref{eq11}:
$$\sum_{k\in I}h_k (f(\x^k)-f(\x^*)) \leqslant\varepsilon\sum_{k\in I}h_k,$$
where for $\hat{\x}^N:=\sum\limits_{k\in I}\frac{h_k\x^k}{\sum_{k\in I}h_k}$ holds~\eqref{th1_1}.
Wherein $\forall k\in I\;g(\x^k)\leqslant\varepsilon||\nabla g(\x^k)||_2\leqslant\varepsilon M_g$ and holds~\eqref{th1_2}.

\end{pf}

Let us estimate the number of iterations necessary to fulfill the stopping criterion \eqref{eq11} in the case of a Lipschitz-continuous objective functional $$|f(x)-f(y)|\leqslant M_f||x-y||_2.$$ It is clear that $\forall k\in I\;||\nabla f(x^k)||_2\leqslant M_f$ and therefore
$$|J|+\sum_{k\in I}\frac{1}{||\nabla f(x^k)||_*^2}\geqslant|J|+\frac{|I|}{M_f^2}\geqslant(|I|+|J|)\frac{1}{\max\{1,M_f^2\}}.$$
This means that for
\begin{eqnarray}
\label{alg2_iter_est}
    \geqslant\frac{2\Theta_0^2\max\{1,M_f^2\}}{\varepsilon^2}
\end{eqnarray}
the stopping criterion \eqref{eq11} is obviously fulfilled, that is, the desired accuracy is achieved in $O\left(\frac{1}{\varepsilon^2}\right)$ iterations.

\subsection{Modification for the logarithm utility functions}

Note that the most common utility functions for networks are logarithms, i.e. $ u_k(x_k) = \log x_k$. However, the logarithm is not a Lipschitz function on $\R^n_{+}$, since its gradient is unlimited near zero. However, consider the following modification of Algorithm~\ref{alg4m}. We shift the boundary of the feasible set from zero, i.e. let $x_k \geq \varepsilon n, \, k = 1, \ldots, n.$ Then, by the definition of the gradient of the logarithm, the utility function will be Lipschitz with the constant $M_U = \frac{1}{\varepsilon}$, i.e.
$$
||\nabla U(\x)||_2 \leq \sum_{k=1}^n | u_k'(x_k)| \leq n \cdot \dfrac{1}{\varepsilon n} = \dfrac{1}{\varepsilon}.
$$
Firstly, to solve this problem, we apply Algorithm~\ref{alg4m} for $N = \left \lceil 2\dfrac{\Theta_0^2}{\varepsilon^4} \right \rceil$ with $h_k=\frac{\varepsilon^2}{\|\nabla f(\x^k)\|_2}$ for $k \in I$ and $h_k=\frac{\varepsilon^2}{\|\nabla g_{j_k}(\x^k)\|_2}$ for $k \not \in I$. Then, we obtain the following estimation for the convergence rate
\begin{cor}
After the $N = \left \lceil 2\dfrac{\Theta_0^2}{\varepsilon^4} \right \rceil$ steps of Algorithm~\ref{alg4m}, the following inequality holds:
$$\min\limits_{k\in I}f(\x^k)-f(\x^*)\leq \varepsilon.$$
\end{cor}

Moreover, let us estimate the convergence rate of Algorithm~\ref{alg1} applied to this problem. Using estimation~\eqref{alg2_iter_est} with $\max$ at $M_{U} = \frac{1}{\varepsilon}$ we obtain the following corollary.
\begin{cor}
After the $N = \left \lceil 2\dfrac{\Theta_0^2}{\varepsilon^4} \right \rceil$ steps of Algorithm~\ref{alg1}, the following inequality holds:
$$f(\hat{\x}^N)-f(\x_*)\leqslant\varepsilon $$
where $\hat{\x}^N=\frac{1}{\sum_{k\in I}h_k}\sum\limits_{k\in I}h_k\x^k$.
\end{cor}

Note that the convergence rates of Algorithm~\ref{alg4m} and Algorithm~\ref{alg1} are of the same order, but due to the adaptability of the stopping criterion, Algorithm~\ref{alg1} works better in practice. Moreover, Algorithm~\ref{alg1} does not require modification of steps.

\section{Ellipsoid method}\label{sect_em}
Consider the transition to the dual problem
for~\eqref{eq_main_task}. Let $ \boldsymbol{\lambda} =
(\lambda_1,\, \ldots,\, \lambda_m) \, \in \, \mathbb{R}^m_+$ be a vector of dual multipliers , which can be interpreted as a compound price vector. Define
dual objective function 
\begin{equation*}
\label{eq_phi_first}
\varphi(\boldsymbol{\lambda}) = \max_{\mathbf{x} \, \in \, \R_+^n} 
\left \{
\sum_{k = 1}^n u_k(x_k) + \langle \boldsymbol{\lambda},
\mathbf{b} - \sum_{k = 1}^n \mathbf{C}_k x_k \rangle
\right \} =
\end{equation*}
\begin{equation*}
= \langle \boldsymbol{\lambda}, \mathbf{b} \rangle + \sum_{k = 1}^n (u_k(x_k(\boldsymbol{\lambda})) - \langle \boldsymbol{\lambda},\, \mathbf{C}_k x_k(\boldsymbol{\lambda}) \rangle ),
\end{equation*}
and users choose the optimal data rates~$x_k$ by solving the following optimization problem
\begin{equation}
\label{eq_xi_argmax}
x_k(\boldsymbol{\lambda}) = \argmax_{x_k \, \in \, \R_+} \left \{
u_k(x_k) - x_k \langle \boldsymbol{\lambda},\, \mathbf{C}_k \rangle
\right \}.
\end{equation}
Then to find the optimal prices $\boldsymbol{\lambda^*}$
we need to solve the problem
\begin{equation}
\label{dual_problem}
\min_{\boldsymbol{\lambda} \, \in \, \R_+^m} \varphi(\boldsymbol{\lambda}).
\end{equation}
Suppose that for the primal problem the Slater condition is satisfied, then due to the strong duality both the primal and the dual problems will have a solution. Using Slater's condition, one can compactify the solution of the dual problem. We assume that the following estimate is correct for solving the dual problem: 
$$
||\boldsymbol{\lambda}^*||_2 \leq R.
$$
In this case, the value $R$ does not affect the operation of the algorithms under consideration, but $R$ is only present in their convergence rate estimations. 

To solve the dual problem, we consider the ellipsoid method.
    \begin{algorithm}
		\caption{Ellipsoid method}\label{alg_ellipsoid}
		\begin{algorithmic}[1]
			\Require $u_k(x_k), k = 1, ..., n$~--- concave functions.
			
			\State $B_0 := 2R \cdot I_n$
			\For{$t=0, ..., N-1$}
    			\State Compute $\nabla \varphi(\boldsymbol{\lambda}^t)$
    			
    			\State $\mathbf{q}_t := B_t^T \nabla \varphi(\boldsymbol{\lambda}^t)$
    			\State $\mathbf{p}_t := \displaystyle \dfrac{B_t^T \mathbf{q}_t}{\sqrt{\mathbf{q}_t^T B_t B_t^T \mathbf{q}_t}}$
    			\State $$B_{t+1} := \displaystyle \frac{m}{\sqrt{m^2 - 1}} B_t + \left(\dfrac{m}{m+1} - \dfrac{m}{\sqrt{m^2 - 1}}\right) B_t \mathbf{p}_t \mathbf{p}_t^T$$
    			\State $\boldsymbol{\lambda}^{t+1} := \boldsymbol{\lambda}^t - \displaystyle \frac{1}{m+1} B_t \mathbf{p}_t$
            \EndFor
			\State\Return $\boldsymbol{\lambda}^N$
		\end{algorithmic}
	\end{algorithm}
As the starting point of the method, we take the zero vector, i.e. $\boldsymbol{\lambda}^0 = 0$. The problem will be solved on the set $\Lambda_{2 R}$, where $$
    \Lambda_{2R} = \{\boldsymbol{\lambda} \in \mathbb{R}_{+}^m: \|\boldsymbol{\lambda}\|_2 \leq 2  R\}.
    $$ 
To restore the solution of the primal problem by the solution of the dual problem, it is necessary to determine the accuracy certificate $\xi $ for the ellipsoid method. The accuracy certificate is a sequence of 
weights $\xi = \{\xi_t\}_{t=0}^{N-1}$ such that
	$$
	\xi_t \geq 0, \sum_{t=0}^{N-1} \xi_t = 1.
	$$
A detailed description of the construction of such a certificate can be found in \cite{nemirovski2010accuracy}.

Now we formulate a theorem of convergence rate estimation~\cite{ivanova2019res}.
\begin{thm}\label{ellips_methods}
Let Algorithm~\ref{alg_ellipsoid} start
with an initial ball $B_0 = \{\boldsymbol{\lambda} \in \mathbb{R}^m: \|\boldsymbol{\lambda}\|_2 \leq 2R\}$. Then after 
	\begin{equation*} \label{th_saga}
	    N = 2m(m+1) \left\lceil\log \left(\frac{32\cdot4 M R}{\varepsilon}\right) \right \rceil
	\end{equation*}
	the following inequalities will hold
	$$
    U(\mathbf{x}^{*}) - U(\hat{\mathbf{x}}^N) \leq \varepsilon, \; \|[C \hat{\mathbf{x}}^N - \mathbf{b}]_{+}\|_2 \leq \varepsilon,
    $$

    where $\displaystyle \hat{\mathbf{x}}^{N} = \sum_{t \in I_{N}} \xi_t \mathbf{x}(\boldsymbol{\lambda}^t), \, \, I_{N} = \left\{ t \leq N-1: \, \boldsymbol{\lambda}^t \in \text{\normalfont int}\;\Lambda_{2R} \right\}$.
\end{thm}

%
%
\section{Experiments}
To test the performance of 
Algorithm~\ref{alg1}
we compared it with the ellipsoid
method (Algorithm~\ref{alg_ellipsoid}).
The
behavior of the methods was tested in problems~\eqref{eq_main_task} of different 
configurations of networks and with different accuracy~$\varepsilon$. 

The routing matrix $C$ was generated
as follows: $C_i^j = 1$ with probability $p=0.5$ or $C_i^j = 0$ with probability $1-p = 0.5$. The elements of the vector~$\mathbf{b}$ are 
uniform random variables: $b_i \, \in \, [0.1, \, 0.4]$. 
The utility functions are logarithmic.
The initial values for Algorithm~\ref{alg1} and the ellipsoid method~(EM) are $\mathbf{x}^0 = 0$ and $\boldsymbol{\lambda}^0 = 10^{-20}$, respectively. The radius~$2R$ of the 
initial ball in the ellipsoid method and the radius~$R = \sqrt{2} \Theta $ of the ball containing $\mathbf{x}^*$ in Algorithm~\ref{alg1} were chosen experimentally in such a way that the intermediate solutions obtained by the methods remained inside the given set at each iteration. The required solution accuracy~$\varepsilon$ was chosen so that the boundary shift~$n \varepsilon$ of the feasible set from zero was small enough, that is no more than $\sim 10^{-1}$.

\setlength\tabcolsep{3pt} 
\begin{table}[htb]
\caption{Convergence results of Algorithm~\ref{alg1}~(A2) and the ellipsoid method~(EM), $\varepsilon = 6e-4$}
\label{tab1}
\begin{tabular}{|l|l|c|c|c|c|c|c|}
\hline
\multirow{2}{*}{}   & \multicolumn{1}{c|}{n} & \multicolumn{2}{c|}{50} & \multicolumn{2}{c|}{100} & \multicolumn{2}{c|}{200} \\ \cline{2-8} 
                    & \multicolumn{1}{c|}{m} & 100        & 150        & 100         & 150        & 100         & 150        \\ \hline
\multirow{2}{*}{A2} & Iter                   & 142243     & 142516     & 171292      & 174270     & 193621      & 198585     \\ \cline{2-8} 
                    & Time, s                & 16.77      & 21.91      & 33.56       & 37.63      & 46.8        & 49.22      \\ \hline
\multirow{2}{*}{EM} & Iter                   & 512749     & 758327     & 531448      & 760537     & 532992      & 761008     \\ \cline{2-8} 
                    & Time, s                & 601.74     & 885.54     & 1022.73     & 1293.15    & 1418.67     & 1481.02    \\ \hline
\end{tabular}
\end{table}

\setlength\tabcolsep{2pt} 
\begin{table}[htb]
\caption{Convergence results of Algorithm~\ref{alg1}~(A2) and the ellipsoid method~(EM), $\varepsilon = 3e-4$}
\label{tab2}
\begin{tabular}{|l|l|c|c|c|c|c|c|}
\hline
\multirow{2}{*}{}   & \multicolumn{1}{c|}{n} & \multicolumn{2}{c|}{50} & \multicolumn{2}{c|}{100} & \multicolumn{2}{c|}{200} \\ \cline{2-8} 
                    & \multicolumn{1}{c|}{m} & 100        & 150        & 100         & 150        & 100         & 150        \\ \hline
\multirow{2}{*}{A2} & Iter                   & 8724510    & 9105234    & 9006192     & 9574296    & 9157003     & 9611472    \\ \cline{2-8} 
                    & Time, s                & 921.38     & 1224.70    & 1276.73     & 1411.67    & 1424.12     & 1670.74    \\ \hline
\multirow{2}{*}{EM} & Iter                   & 603578     & 801775     & 628267      & 833323     & 633571      & 850051     \\ \cline{2-8} 
                    & Time, s                & 1084.22    & 1317.61    & 1321.48     & 1705.46    & 1492.88     & 1921.06    \\ \hline
\end{tabular}
\end{table}

\setlength\tabcolsep{6pt}
\begin{table}[htb]
\caption{Convergence results of Algorithm~\ref{alg1}~(A2) and the ellipsoid method~(EM), $\varepsilon = 2e-4$}
\label{tab3}
\begin{tabular}{|l|l|c|c|ll}
\hline
\multirow{2}{*}{}                         & \multicolumn{1}{c|}{n} & \multicolumn{2}{c|}{50}  & \multicolumn{2}{c|}{100}                                      \\ \cline{2-6} 
                                          & \multicolumn{1}{c|}{m} & 100         & 150        & \multicolumn{1}{c|}{100}      & \multicolumn{1}{c|}{150}      \\ \hline
\multirow{2}{*}{A2}                       & Iter                   & 25225735    & 29752323   & \multicolumn{1}{c|}{26055762} & \multicolumn{1}{c|}{33846145} \\ \cline{2-6} 
                                          & Time, s                & 1367.54     & 1569.25    & \multicolumn{1}{c|}{1550.62}  & \multicolumn{1}{c|}{1796.34}  \\ \hline
\multirow{2}{*}{EM}                       & Iter                   & 599423      & 960529     & \multicolumn{1}{c|}{618783}   & \multicolumn{1}{c|}{971525}   \\ \cline{2-6} 
                                          & Time, s                & 1223.86     & 1515.32    & \multicolumn{1}{c|}{1677.25}  & \multicolumn{1}{c|}{1900.03}  \\ \hline
\multicolumn{1}{|c|}{\multirow{2}{*}{}}   & \multicolumn{1}{c|}{n} & \multicolumn{2}{c|}{200} &                               &                               \\ \cline{2-4}
\multicolumn{1}{|c|}{}                    & \multicolumn{1}{c|}{m} & 100         & 150        &                               &                               \\ \cline{1-4}
\multicolumn{1}{|c|}{\multirow{2}{*}{A2}} & Iter                   & 28532359    & 37244837   &                               &                               \\ \cline{2-4}
\multicolumn{1}{|c|}{}                    & Time, s                & 1723.07     & 1985.52    &                               &                               \\ \cline{1-4}
\multicolumn{1}{|c|}{\multirow{2}{*}{EM}} & Iter                   & 667294      & 1021528    &                               &                               \\ \cline{2-4}
\multicolumn{1}{|c|}{}                    & Time, s                & 1891.18     & 2293.34    &                               &                               \\ \cline{1-4}
\end{tabular}
\end{table}

The results of the experiments are presented in Tables~\ref{tab1}-\ref{tab3}.
As one can see from Tables~\ref{tab1}-\ref{tab2}, for $\varepsilon = 6e-4$ and $\varepsilon = 3e-4$, the proposed algorithm shows better time than the ellipsoid method. Even for
high solution accuracy, $\varepsilon = 2e-4$,
Algorithm~\ref{alg1} showed 
a large number of iterations and
almost the same time as EM, as shown in Table~\ref{tab3}.
So, in a case where very high solution accuracy is not required, it is reasonable to apply the proposed
algorithm.

The conducted experiments confirm the following 
theoretical
fact about Algorithm~\ref{alg1}:
the convergence rate of Algorithm~\ref{alg1} depends only on 
the smoothness level of the target function
and of the constraints and does
not depend on the number of constraints~$m$
(see Section~\ref{sect_mir}).
Unlike in the case of the ellipsoid method, where there is a quadratic growth over $m$
in the theoretical number of iterations
(Theorem~\ref{ellips_methods}).
So, if we compare the number of iterations ({\it Iter})
for the same $n$ and for $m = 100$ and $m = 150$
in Tables~\ref{tab1}-\ref{tab3},
one can notice that the number of iterations for EM
increases almost $1.5-2$ times and
it is not the same for Algorithm~\ref{alg1}. 
The theoretical results tell us that for Algorithm~\ref{alg1}, as $m$ increases and $n$ is the same, the number of iterations should not change. But due to the adaptability of the stopping criterion, in practice it changes slightly, since in both cases this number is less than the theoretical convergence rate estimate.

%
%
\section{Conclusion}

In conclusion, we note that despite
the theoretical attractiveness of Algorithm~\ref{alg4m} (for example, as can be seen from Section~\ref{Mirr_Des_1}, one can obtain estimates for 
cases with an inexact oracle and objective
functions
with different smoothness levels), our
experiments showed that
Algorithm~\ref{alg1} is significantly faster (both in time and in number of iterations) in practice due to adaptability of the stopping criterion. 

Moreover, we considered Algorithm~$3$ from~\cite{sasb_adap_mir2018} (see also Algorithm $1$ from \cite{baynd_lec_notes2018}). Note, in comparison with Algorithm~\ref{alg1} this algorithm guarantees a better estimate for the residual by constraint $g(\x) \leq \varepsilon$ with similar estimates for the objective function and similar complexity $O \left(\varepsilon^{-2}\right)$. 
However, in practice, it works much worse than  Algorithm~\ref{alg1} and Algorithm~\ref{alg4m}.

Another important observation concerning the comparison of Algorithms~\ref{alg1} and \ref{alg4m} with methods from~\cite{baynd_lec_notes2018, sasb_adap_mir2018} in practice is that for medium and large networks. 
The practical result in a reasonable time can be obtained only from Algorithm~\ref{alg1}.

In conclusion we note that  approaches considered by us allow us to determine the prices of connections. Here you can use the properties of primal-duality of considered  Mirror Descent methods
\cite{baynd_zhvm2018, baynd_lec_notes2018}.

\begin{ack}
The authors are very grateful to P. Dvurechensky and Yu.~Nesterov for fruitful discussions.
\end{ack}

\bibliography{our}             
                                                   







\end{document}